\magnification=\magstep1
\pageno=1
\input amstex
\documentstyle{amsppt}
\hsize 6truein
\nologo
\pagewidth{32pc}
\pageheight{45pc}
\mag=1200
\baselineskip=15 pt
\topmatter
\title
Finite Bruck loops\endtitle
\author Michael Aschbacher, Michael K. Kinyon, and J.D. Phillips\endauthor
\affil California Institute of Technology, Indiana University South Bend,\\ and Wabash College \endaffil
\address Pasadena, California 91125\endaddress
\address South Bend, Indiana 46634\endaddress
\address Crawfordsville, Indiana 47933\endaddress
\thanks The first author was partially supported by NSF-0203417\endthanks
\endtopmatter

\document

Let $X$ be a magma; that is $X$ is a set together with a binary operation $\circ$ on $X$.
For each $x\in X$ we obtain maps $R(x)$ and $L(x)$ on $X$ defined by
$R(x) : y\mapsto y\circ x$ and $L(x) : y\mapsto x\circ y$ called {\it right and left translation}
by $x$, respectively. A {\it loop} is a magma $X$ with an identity 1 such that $R(x)$ and $L(x)$
are permutations of $X$ for all $x\in X$. In essence loops are groups without the associative axiom.
See \cite{Br, Pf} for further discussion of basic properties of loops.

Certain classes of loops have received special attention:
 A loop $X$ is a (right) {\it Bol loop} if it satisfies the (right) {\it Bol identity} (Bol):
$$
((z\circ x)\circ y)\circ x = z\circ ((x\circ y)\circ x). \tag{Bol}
$$
or equivalently
$$
R(x)R(y)R(x) = R((x\circ y)\circ x). \tag{Bol2}
$$
for all $x,y,z\in X$.
In a Bol loop, the subloop $\langle x\rangle$ generated by $x\in X$ is a group.
Thus we can define $x^{-1}$ and the order $|x|$ of $x$ to be, respectively,
the inverse of $x$ and the order of $x$ in that group. For basic facts about Bol
loops, see \cite{Ro}.

A loop $X$ in which inverses are defined satisfies the {\it automorphic inverse property (AIP)}
if $(x\circ y)^{-1}=x^{-1}\circ y^{-1}$ for all $x,y\in X$. Finally $X$ is a {\it Bruck loop}
if $X$ is a Bol loop satisfying the AIP. Bruck loops are also known as K-loops \cite{Ki} and
gyrocommutative gyrogroups (see, e.g., \cite{FU}).

We prove many of our results on loops by translating them into results about groups, using an
observation of Reinhold Baer in  \cite{Ba}:
Given a loop $X$, define $K=\{R(x) : x\in X\}$ regarded as a subset of the
symmetric group $Sym(X)$ on $X$, $G=\langle K\rangle$ to be the subgroup
of $Sym(X)$ generated by $K$, and $H=G_1$ to be the stabilizer in $G$ of the identity 1 of $X$.
Set $\epsilon(X)=(G,H,K)$.
 In the loop theory literature, $G$ is usually called the {\it right multiplication group}
of $X$, but since we will make no reference to groups with left translations as generators,
we will follow \cite{A2} and simply call $G$ the {\it enveloping group} of $X$. The
subgroup $H$ is the (right) {\it inner mapping group} of $X$, and we call $\epsilon(X)$
the {\it envelope} of $X$.

We can now state our main theorems:
\smallpagebreak

\proclaim{Theorem 1}
Let $X$ be a finite Bruck loop with enveloping group $G$. Then

(1) $X=O^{2'}(X)*O(X)$ and $G=O^{2'}(G)*O(G)$ are central products.

(2) $O^{2'}(X)\cap O(X)=Z$ is a subloop of $Z(X)$ of odd order and
      $O^{2'}(G)\cap O(G)$ is a subgroup of $Z(G)$ of odd order.

(3) $X/Z=O^{2'}(X)/Z\times O(X)/Z$.

(4) $O^2(X)=O(X)$, so $O^{2'}(X)/Z$ is a 2-element loop.

(5) If $X$ is solvable then $O^{2'}(X)=O_2(X)$, so
      $X=O_2(X)\times O(X)$ and $G=O_2(G)\times O(G)$.
\endproclaim
\smallpagebreak

See \cite{FGT} for notation and terminology involving groups.
An element of $X$ of order a power of $2$ is called a $2${\it -element}, and
we say that $X$ is a $2${\it -element loop} if every element of $X$ is a $2$-element.
Given a set $\pi$ of primes, $X$ is a $\pi$-{\it loop} if $\pi(|X|)\subseteq\pi$.
In particular, if $|X|$ is a power of $2$, then $X$ is a $2${\it -loop}, while if $|X|$ is
odd, then $X$ is a $2'${\it -loop}. We write $O_{\pi}(X)$
for the largest normal $\pi$-subloop of $X$, so that $O_2(X)$ is the largest normal
$2$-subloop. Further we abbreviate $O_{2'}(X)$, the largest normal subloop
of odd order, by $O(X)$. Finally we write $O^2(X)$, $O^{2'}(X)$ for the subloop
generated by all elements of odd order, all $2$-elements, respectively.

The {\it center} $Z(X)$ of a loop $X$ is the set of all $a\in X$ such that
$a\circ (x\circ y) = x\circ (a\circ y) = (x\circ a)\circ y = x\circ (y\circ a)$ for
all $x,y\in X$. The center is a normal subloop (see Section 1 for the definition
of normality).

One of the main tools in the proof of Theorem 1 is the following result about
arbitrary Bruck loops, which is of independent interest:
\smallpagebreak

\proclaim{Theorem 2}
Let $X$ be a Bruck loop and $x,y\in X$ with
$x$ a 2-element and $y$ an element of odd order.
Then $R(x)R(y)=R(x\circ y)=R(y\circ x)=R(y)R(x)$.
Hence $x\circ y=y\circ x$.
\endproclaim
\smallpagebreak

Theorem 1 reduces the study of finite Bruck loops to the study of 2-element loops and loops of odd order.
 The category of Bruck loops of odd order is essentially equivalent to the category of pairs $(G,\tau )$ where
$G$ is a group of odd order and $\tau$ an involutory automorphism of $G$ such that $G=[G,\tau ]$ and
$C_{Z(G)}(\tau)=1$. This fact goes back to Glauberman in \cite{G2} and \cite{G3}; see also 5.7, 5.8, and 5.10.
 As a result, Bruck loops of odd order are well behaved and well understood.
 On the other hand Bruck 2-loops and Bruck 2-element loops seem difficult to analyze.

It seems possible that all finite Bruck loops $X$ are solvable, and hence $X=O_2(X)\times O(X)$.
 Our next theorem is a step toward proving that finite Bruck loops are indeed solvable.
 Define a finite Bruck loop to be an {\it M-loop} if each proper section of $X$ is solvable,
 but $X$ is not solvable.
\smallpagebreak

\proclaim{Theorem 3} Let $X$ be an M-loop, $\epsilon(X)=(G,H,K)$, $J=O_2(G)$, and $G^*=G/J$.
 Then

(1) $X$ is a simple $2$-element loop.

(2) $G^*\cong PGL_2(q)$, with $q=2^n+1\geq 5$, $H^*$ is a Borel subgroup of $G^*$, and
$K^*$ consists of the involutions in $G^*-F^*(G^*)$.

(3) $F^*(G)=J$.

(4) Let $n_0=|K\cap J|$ and $n_1=|K\cap aJ|$ for $a\in K-J$.
 Then $n_0$ is a power of $2$, $n_0=n_12^{n-1}$, and $|X|=|K|=(q+1)n_0=n_12^n(2^{n-1}+1)$.
\endproclaim
\smallpagebreak

One would like to show that M-loops do not exist, and hence show that finite Bruck loops are solvable.
 Theorem 3 identifies a set of obstructions to that goal.
 This is essentially the same set of obstructions to the Main Theorem of \cite{A2} on Bol loops of exponent $2$.
 See section 12 of \cite{A2} for a discussion of possible approaches to eliminating these obstructions or
 alternatively to constructing examples of M-loops. These approaches involve the analysis of Bruck $2$-loops.

A loop $X$ is said to be an $A_r$-{\it loop} if its inner mapping group acts as a group
of automorphisms of $X$ in its permutation representation on $X$.
The class of finite Bol loops which are also $A_r$-loops is much larger than the class of
finite Bruck loops; for example the former class includes all finite groups.
 Still (cf. Lemma 8.1) the latter class can be described in terms of the former
class and the class of finite groups, allowing us to prove:
\smallpagebreak

\proclaim{Corollary 4} Let $X$ be a finite loop which is both a Bol loop and an $A_r$-loop.
 Then $X$ is solvable iff the enveloping group of $X$ is solvable.
\endproclaim
\smallpagebreak

The proof of Theorem 1 uses the solvability of groups of odd order, established by Feit and Thompson
 in \cite{FT}, Glauberman's $Z^*$-Theorem \cite{G1}, and several other results from the theory of finite
 groups, whose proofs are a bit easier and can be found in \cite{FGT}.
 The proof of Theorem 3 involves appeals to the Main Theorem of \cite{A2} and its proof, which in turn
 uses the classification of the finite simple groups, together with deep knowledge of the subgroup structure
 of the automorphism groups of those groups.
\smallpagebreak

\subheading{Section 1}
 Loops, folders, envelopes, and twisted subgroups
\smallpagebreak

In this section we recall some notation and terminology involving loops, summarize some facts about loops,
and references for those facts.

In \cite{A2}, a {\it loop folder} is defined to be a triple $\xi = (G,H,K)$, where $G$ is a group,
$H$ is a subgroup of $G$, $K$ is a subset of $G$ containing 1, and for all $g\in G$, $K$ is a set
of coset representatives for $H^g$ in $G$. The folder is an {\it envelope} if $G=\langle K\rangle$
and {\it faithful} if $\ker_H(G)=1$, where $\ker_H(G)$ is the largest normal subgroup of $G$
contained in $H$.

For example if $X$ is a loop then $\epsilon(X)$ is a faithful loop envelope.

Section 1 of \cite{A2} contains the definition of a category of loop folders and functors
$\epsilon$ and $l$ to and from the category of  loops and the category of loop folders.
 The reader is directed to \cite{A2} for notation, terminology, and results about folders and these functors.

A {\it twisted subgroup} of a group $G$ is a subset $K$ of $G$ such that $1\in K$ and for all
$x,y\in K$, $xy^{-1}x\in K$. See (\cite{A2}, \S5) for a brief discussion of twisted subgroups
taken from \cite{A1}.

A folder $\xi = (G,H,K)$ is a {\it Bol loop folder} if $K$ is a twisted subgroup of $G$.
 Further (cf. \cite{A2}, 6.1) a loop $X$ is a Bol loop iff $\epsilon(X)=(G,H,K)$ is a Bol loop folder.
 In that event there is a normal subgroup $\Xi_K(G)$ of $G$ called the {\it K-radical} of $G$, and a
corresponding normal subloop $\Xi(X)$ of $X$ (which is a group) called the {\it radical} of $X$.
 Moreover if $\Xi_K(G)=1$ then there is a unique automorphism $\tau = \tau_X$ of $G$ such that
$\tau^2=1$ and $K\subseteq K(\tau)$, where
$$
K(\tau )=\{g\in G : g^{\tau} = g^{-1}\}.
$$
See (\cite{A2}, \S6) for further discussion.

Next $X$ is an $A_r$-loop iff $H$ acts on $K$ via conjugation (cf. 4.1 in \cite{A2}).
 Further $X$ is a Bruck loop iff $X$ is a radical free (ie. $\Xi(X)=1$) $A_r$-loop.
 (cf. 6.7 in \cite{A2})

The material in the remainder of this section is elementary and easy.
See, for example, (\cite{Br}, IV.1) for further discussion and proofs.

A {\it normal subloop} of a loop $X$ is the kernel of a loop homomorphism.
 Further a subloop $Y$ of $X$ is normal iff for all $a,b\in X$,
$$a\circ (Y\circ b)=Y\circ (a\circ b) = (a\circ Y)\circ b,$$
in which case the cosets $Y\circ x$, $x\in X$, form the equivalences classes
of an equivalence relation (congruence) on $X$, and we can form the factor loop
$X/Y$ on this set of equivalence classes, with multiplication defined by
$$
(X\circ a)\circ (X\circ b)=X\circ (a\circ b).
$$
Also we obtain the surjective loop homomorphism $\pi : X\to X/Y$ with $x\pi = Y\circ x$ and $\ker(\pi)=Y$.
 We have the usual facts:
\smallpagebreak

\proclaim\nofrills{(1.1)\usualspace } If $\varphi : X\to X'$ is a surjective loop homomorphism with $\ker(\varphi )=Y$ then

(1) $\psi : X/Y\to X'$ defined by $(Y\circ x)\psi = x\varphi$ is an isomorphism with $\pi\psi = \varphi$.

(2) If $U\ \trianglelefteq\ X$ then $U\varphi\ \trianglelefteq\ X'$.

(3) If $V\ \trianglelefteq\ X'$ then $V\varphi^{-1}\ \trianglelefteq\ X$.
\endproclaim

\smallpagebreak

\subheading{Section 2}
 Normal structure of loops
\smallpagebreak

In this section $\xi = (G,H,K)$ is a loop envelope and $X=l(\xi)$.
\smallpagebreak

\proclaim\nofrills{(2.1)\usualspace } Let $\xi_i=(G_i,H_i,K_i)$ be normal subfolders of $\xi$, set $\xi^i=\xi/\xi_i$,
and let $\pi_i : \xi\to\xi^i$ be the natural map of (\cite{A2}, 2.6) with $\ker(\pi_i)=\xi_i$.
 Then

(1) $\xi_{3-i}\pi_i\ \trianglelefteq\ \xi^i$.

(2) Let $G_3=G_1G_2$, $H_3=H\cap G_3$, and $K_3=K\cap G_3$.
 Then $\xi_3=(G_3,H_3,K_3)\ \trianglelefteq\ \xi$ and $\xi_3\pi_i=\xi_{3-i}\pi_i$.

(3) $\xi/\xi_3\cong\xi^i/\xi_{3-i}\pi_i$.

(4) Let $G_0=G_1\cap G_2$, $H_0=H_1\cap H_2$, and $K_0=K_1\cap K_2$.
 Then $\xi_0\ \trianglelefteq\ \xi$.

(5) Set $\bar G=G/G_0$.
 Then $\bar G_3=\bar G_1\times\bar G_2$, $\bar H_3=\bar H_1\times\bar H_2$, and $\bar K_3=\bar K_1\times\bar K_2$.

(6) Let $X_i=l(\xi_i)$.
 Then $X_i$ is normal in $X$ for each $i$ and $X_3/X_0\cong X_1/X_0\times X_2/X_0$.
\endproclaim
\demo{Proof} Let $G^*=G\pi_1$.
 Then $\xi^1=(G^*,H^*,K^*)$ and $\xi_2\pi_1=\xi_2^*=(G_2^*,H_2^*,K_2^*)$.
 By (\cite{A2}, 2.9), $X_i\ \trianglelefteq\ X$ and $\psi_i=l(\pi_i) : X\to X^i=X/X_i$ is an isomorphism.
 As $X_2\ \trianglelefteq X$, $X_2\psi_1\ \trianglelefteq\ X^1$ by 1.1, and then by another application of 1.1,
the preimage $Y$ of $X_2\psi_1$ in $X$ under $\psi_1$ is also normal in $X$.
 By (\cite{A2}, 2.9.3), there is a normal subfolder $\mu = (G_Y,H_Y,K_Y)$ of $\xi$ with $l(\mu)=Y$.
 As $Y\psi_1=X_2\psi_1$, $K_Y^*=K_2^*$, so $s_{\xi^*}(K_2^*)=s_{\xi^*}(K_Y^*)$ is a subenvelope
of $\xi^*$ and hence as $G_2^*=H_2^*K_2^*$, $\xi_2^*$ is a subfolder of $\xi^*$ by (\cite{A2}, 2.1).

Let $k_2\in K$, $k\in K$, and $g\in G$.
 As $\xi_2\ \trianglelefteq\ \xi$ the normality condition (NC) from (\cite{A2}, \S2) is satisfied by $\xi_2$,
so there is $l\in H^g\cap G_2$ and $k_3\in K$ with $k_2k=lk_3$.
 Then $k_2^*k^*=l^*k_3^*$, so $\xi_2^*$ satisfies (NC) in $\xi^*$, establishing (1).

Let $\alpha : \xi^*\to \xi^*/\xi_2^*$ be the natural map and $\beta = \pi\alpha$.
 Then $\beta :\xi\to \xi^*/\xi_2^*$ is a surjective morphism with kernel $\xi_3$, so
$\xi_3$ is a normal subfolder of $\xi$, establishing (2) and (3).

Let $u,g\in G_0$. Then $u=h_uk_u$ with $h_u\in H^g$ and $k_u\in K$.
 As $\xi_i$ is a subfolder of $\xi$ for $i=1,2$, $h_u\in H_i^g$ and $k_u\in K_i$, so
$h_u\in H_1^g\cap H_2^g=H_0^g$ and $k_u\in K_1\cap K_2=K_0$, and hence $\xi_0$ is a subfolder of $\xi$.
 Similarly if $k_0\in K_0$ and $k\in K$ then $k_0k=lk_4$ for some $l\in H^g$ and $k_4\in K$, and as
$\xi_i\ \trianglelefteq\ \xi$ for $i=1,2$, $l\in H_i^g$, so $l\in H_0^g$ and hence $\xi_0$ satisfies (NC) in $\xi$.
 This establishes (4).

Of course
$$
\bar G_3=\bar G_2\times \bar G_2=\bar H_1\bar K_1\times\bar H_2\bar K_2 =
(\bar H_1\times\bar H_2)(\bar K_1\times\bar K_2),
$$
 with $H_1H_2\leq H_3$. Let $a_i\in K_i$, $i=1,2$.
 Then $a_1a_2=hk$, $h\in H_3$, $k\in K_3$.
 Also $a_2^*=a_1^*a_2^*=h^*k^*$, so $h^*=1$; that is $h\in H_1$.
 By symmetry $h\in H_2$, so $\bar h=1$ and hence $\bar K_1\times\bar K_2\subseteq\bar K_3$.
 Therefore $\bar H_3=\bar H_1\times\bar H_3$ and $\bar K_3=\bar K_1\times\bar K_2$, establishing (5).

Finally, applying the functor $l$ to (5), we obtain (6).
\enddemo
\smallpagebreak

\proclaim\nofrills{(2.2)\usualspace } Assume $X$ is finite and $\pi$ is a set of primes.
 Then

(1) $X$ has a largest normal $\pi$-subloop $O_{\pi}(X)$.

(2) There is a normal subfolder $\xi_{\pi}=(G_{\pi},H_{\pi},K_{\pi})$ of $\xi$ such that $l(\xi_{\pi})=O_{\pi}(X)$.
\endproclaim
\demo{Proof} Part (2) is a consequence of (1) and (\cite{A2}, 2.9.3).
 By 2.1.6, if $X_1$ and $X_2$ are normal $\pi$-subloops of $X$ then there is a normal
$\pi$-subloop $X_3$ of $X$ containing $X_1$ and $X_2$, so (1) holds.
\enddemo
\smallpagebreak

\subheading{Section 3}
 Radical free Bol loops.
\smallpagebreak

In this section $X$ is a radical free Bol loop and $\xi = \epsilon(X)=(G,H,K)$.
 Adopt Notation 6.3 from \cite{A2}, and assume:
\smallpagebreak

 \flushpar {\bf Hypothesis 3.1.} $M^+$ is a maximal subgroup of $G^+$ containing
$H\langle\tau\rangle$ and $J^+$ is a normal subgroup of $G^+$ contained in $M^+$.
 Set $G^{+*}=G^+/J^+$, $K_M=K\cap M^+$, and $\Lambda_M=\tau K_M$.
 Let $\Delta = G^+/M^+$ and represent $G^+$ on $\Delta$ via right multiplication.
\smallpagebreak

\proclaim\nofrills{(3.2)\usualspace } (1) $|K|=|\Lambda |$ and $|K_M|=|\Lambda_M|$.

(2) $\Lambda_M=\Lambda\cap M^+$.

(3) $|\Lambda |=|G:M||\Lambda_M|$.

(4) Let $\lambda_i^*$, $1\leq i\leq r$ be representatives for the orbits of $G^*$ on $\Lambda^*$,
$m_i=|\lambda_i^{*G^*}|$, $n_i=|\Lambda\cap\lambda_iJ^+|$, and $n_0=|\Lambda\cap J^+|$.
 Then
$$|\Lambda |=n_0+\sum_{i=1}^rn_im_i.$$

(5) If $|G^+:M^+|$ is odd then $n_0=0$ and each member of $\Lambda$ fixes a unique point of $\Delta$.
\endproclaim
\demo{Proof} As the map $k\mapsto \tau k$ is a bijection of $K$ with $\Lambda$ and $\tau\in M^+$, (1) and (2) hold.
 The proof of (3) is straightforward and is the same as that of 12.5.1 in \cite{A2}.
 Similarly the proofs of (4) and (5) are essentially the same as those of parts (2) and (3) of (\cite{A2}, 12.5),
but we repeat the details for completeness:

Let $\Lambda_i=\{\lambda\in\Lambda : \lambda^*\in\lambda_i^{*G}\}$ and $\Lambda_0=\Lambda\cap J^+$.
 Then $\{\Lambda_i : 0\leq i\leq r\}$ is a partition of $\Lambda$ with $|\Lambda_0|=n_0$ and $|\Lambda_i|=n_im_i$
for $1\leq i\leq r$.
 Thus (4) holds.

Finally assume $\alpha = |G^+:M^+|$ is odd. Then by Sylow's Theorem we may choose $\lambda_i\in M^+$.
 Set
$$
t_i=|\lambda_i^{*G}\cap M^{+*}|.
$$
Then arguing as in the proof of (4),
$$
|\Lambda_M|=n_0+\sum_{i=1}^rn_it_i.\tag{*}
$$
Therefore by (3), (4), and (*),
$$
n_0+\sum_{i=1}^rn_im_i=|\Lambda |=\alpha |\Lambda_M|=\alpha n_0+\sum_{i=1}^r\alpha n_it_i,
$$
so
$$
(\alpha - 1)n_0+\sum_{i=1}^rn_i(\alpha t_i-m_i)=0.\tag{**}
$$
We next claim:
\smallpagebreak

\flushpar (!) For each $1\leq i\leq r$, $\alpha t_i\geq m_i$, with equality iff each $\lambda\in\Lambda_i$
fixes a unique point of $\Delta$.
\smallpagebreak

\flushpar Namely each $\lambda_i\in\Lambda_i^*$ is in some conjugate $M^{+*g}$ of $M^{+*}$, and $|M^{+*g}\cap\Lambda_i|=t_i$, so as
$$
|M^{+*G}|=|G^+:N_G(M^+)|=|G:M^+|=\alpha ,
$$
$m_i=|\Lambda_i^*|\leq\alpha t_i$ with equality iff each $\lambda^*$ is contained in a unique
conjugate of $M^{+*}$ iff $\lambda$ fixes a unique point of $\Delta$.

Finally by (!):
\smallpagebreak

\flushpar (!!) $\sum_{i=1}^rn_i(\alpha t_i-m_i)\geq 0$ with equality iff each $\lambda\in\Lambda-\Lambda_0$
fixes a unique point of $\Delta$.
\smallpagebreak

As $\alpha > 1$ we conclude from (**) and (!!) that (5) holds.
\enddemo
\smallpagebreak

\proclaim\nofrills{(3.3)\usualspace } Adopt the notation of 3.2, assume $1\leq n_0$ is a power of $2$, and
$p$ is an odd prime such that $m_i\equiv 0\mod p$ for all $1\leq i\leq r$.
 Then $H$ contains a Sylow $p$-subgroup of $G$.
\endproclaim
\demo{Proof} By 3.2.4, $|\Lambda |\equiv n_0\mod p$, so as $1\leq n_0$ is a power of $2$
and $p$ is odd, $|\Lambda |$ is relatively prime to $p$. Thus as $|G:H|=|\Lambda |$, the lemma follows.
\enddemo
\smallpagebreak

\subheading{Section 4}
 The proof of Theorem 2
\smallpagebreak

In this section $X$ is a Bruck loop and $u,v\in X$.
\smallpagebreak

\proclaim\nofrills{(4.1)\usualspace } (1) The map $x\mapsto R(x)$ is an isomorphism
of $\langle v\rangle$ with $\langle R(v)\rangle$.

(2) For all integers $m, n$, $(u\circ v^m)\circ v^n=u\circ v^{m+n}$.
\endproclaim
\demo{Proof} Part (1) is well known (cf. \cite{A2}, 6.8), and (2) is just a restatement of (1).
\enddemo
\smallpagebreak

\proclaim\nofrills{(4.2)\usualspace } $(u\circ v)^2=(v\circ u^2)\circ v$.
\endproclaim
\demo{Proof} This appears in Lemma 1 in \cite{G2}, but we supply a proof for completeness.
Let $w\in X$ and set $z = w^{-1}\circ x^{-1}$. Since $X$ has the AIP, using 4.1.2:
$$
z\circ (w\circ x)^2 = (w\circ x)^{-1}\circ (w\circ x)^2 = w\circ x .
$$
Next
$$
w\circ x = w^{-1} R(x^{-1})R(x) R(w^2) R(x) = z R(x)R(w^2)R(x) = z R((x\circ w^2)\circ x),
$$
using 4.1 and the Bol identity (Bol2). Thus $z \circ (w\circ x)^2 = z\circ ((x\circ w^2)\circ x)$,
and cancelling $z$, we obtain the lemma.
\enddemo
\smallpagebreak

\proclaim\nofrills{(4.3)\usualspace } Assume $k$ is a positive integer such that $u^{2^j}$ commutes
with $v^{2^{j-1}}$ for each $1\leq j\leq k$.
 Then $(u\circ v)^{2^i}=u^{2^i}\circ v^{2^i}$ for each $0\leq i\leq k$.
\endproclaim
\demo{Proof} The lemma is trivial if $i=0$. When $i=1$, 4.1.2 and 4.2 say
$$(u\circ v)^2=(v\circ u^2)\circ v=(u^2\circ v)\circ v=u^2\circ v^2.$$
Finally complete the proof by induction on $i$, using the validity of the lemma at $i=1$.
\enddemo
\smallpagebreak

\proclaim\nofrills{(4.4)\usualspace } If $v\in\langle u\circ v\rangle$ then $R(u)R(v)=R(u\circ v)=R(v\circ u)=R(v)R(u)$.
\endproclaim
\demo{Proof} As $v\in\langle u\circ v\rangle$, $u=(u\circ v)\circ v^{-1}\in\langle u\circ v\rangle$ by 4.1.2.
 Then the lemma follows from 4.1.1.
\enddemo
\bigpagebreak

With these lemmas in hand, we can prove Theorem 2. Let $x,y\in X$ with $|x|=2^n$ and $|y|$ odd.
 We prove
$$
R(x)R(y)=R(x\circ y)=R(y\circ x)=R(y)R(x)\tag{*}
$$
by induction on $n$. Observe (*) implies $x\circ y=y\circ x$ as $R$ is injective.

When $n=0$, (*) is trivial. Assume $n>0$ and (*) holds for $i<n$.
 Then as $|x^2|=2^{n-1}$, each element of $\langle x^2\rangle$ commutes
with each element of $\langle y\rangle$ by the induction assumption.
 Therefore by 4.3, $(x\circ y)^{2^n}=x^{2^n}\circ y^{2^n}=y^{2^n}$,
 so $y^{2^n}\in\langle x\circ y\rangle$. Then as $|y|$ is odd, $y\in\langle x\circ y\rangle$,
so (*) holds by 4.4.

This completes the proof of Theorem 2.
\smallpagebreak

\subheading{Section 5}
 Bruck loops
\smallpagebreak

\proclaim\nofrills{(5.1)\usualspace } Let $X$ be a loop with envelope $\xi = (G,H,K)$.
 Then the following are equivalent:

(1) $X$ is a Bruck loop.

(2) $X$ is an $A_r$-loop and $X$ is radical free.

(3) $H$ acts via conjugation on $K$ and $\Xi_K(G)=1$.

(4) $\Xi_K(G)=1$ and $H\leq C_G(\tau_X)$.
\endproclaim
\demo{Proof} Parts (1) and (2) are equivalent by (\cite{A2}, 6.6).
 Assume (2). Then $X$ is radical free, so $\Xi_K(G)=1$ by definition.
 As $X$ is an $A_r$-loop, $H$ acts on $K$ by 4.1 in \cite{A2}, so (3) holds.
 The proof of 6.7 in \cite{A2} shows that (3) implies (4).
 Finally the proof of (4) in 6.6 of \cite{A2} shows that (4) implies
$X$ is an $A_r$-loop; thus (4) implies (2).
\enddemo
\smallpagebreak

A loop folder $\xi = (G,H,K)$ is a {\it Bruck loop folder} if $\xi$ is a Bol loop folder,
$\Xi_K(\langle K\rangle )=1$, and $H$ acts on $K$ via conjugation.

In the remainder of the section assume $\xi =(G,H,K)$ is a finite Bruck loop folder.
 We adopt the following notational conventions:
\smallpagebreak

\flushpar {\bf Notation 5.2.} As $\Xi_K(\langle K\rangle)=1$, from (\cite{A2}, 5.1.3c), there
is a unique automorphism $\tau = \tau_{\xi }$ of $\langle K\rangle$ such that $\tau^2=1$ and $K\subseteq K(\tau)$.
 As $H$ acts on $K$, $H\cap\langle K\rangle$ centralizes $\tau$ by the uniqueness of $\tau$.
 As $\xi$ is a loop folder, $K$ is a set of coset representatives for $H$ in $G$, so as $\tau$
centralizes $H\cap\langle K\rangle$ there is a unique extension of $\tau$ to $G$ defined by
$\tau : hk\mapsto hk^{\tau}$ for $h\in H$ and $k\in K$.
Form the semidirect product $G^+=G\langle \tau\rangle$ of $G$ by $\tau$ and let
$\Lambda = \tau K\subseteq G^+$. By (\cite{A2}, 5.1), $\Lambda$ is $\langle K\rangle$-invariant,
so as $H$ centralizes $\tau$ and acts on $K$, and as $G=HK$, $\Lambda$ is also $G$-invariant.

Let $G^+_{\tau}=C_{G^+}(\tau)$, $G_{\tau}=C_G(\tau)$, and $K_{\tau}=C_K(\tau)$.
 Set $\xi_{\tau}=(G_{\tau},H,K_{\tau})$.

For $U\subseteq H$, let $K_U=C_K(U)$, $G_U=N_G(U)$, $H_U=N_H(U)$, and $\xi_U=(G_U,H_U,K_U)$.
\smallpagebreak

\proclaim\nofrills{(5.3)\usualspace } (1) $\xi$ is an $A_r$-loop folder; that is $H$ acts on $K$ via conjugation.

(2) For each $k\in K$, $H\cap H^k=C_H(k)$.

(3) $H$ controls $G$-fusion in $H$.

(4) If $\mu=(G_{\mu},H_{\mu},K_{\mu})$ is a subfolder of $\xi$ then $\mu$ is a Bruck loop folder,
$\tau$ acts on $G_{\mu}$, and $\tau_{\mu}=\tau_{|G_{\mu}}$.

(5) Suppose $\pi : \xi\to\eta = (G_{\eta},H_{\eta},K_{\eta})$ is a surjective homomorphism of loop folders,
and let $\xi_0=(G_0,H_0,K_0)=\ker(\pi)$. Then $\tau$ acts on $G_0$, $\eta$ is a Bruck loop folder,
$\tau_{\eta}=\tau '$, where $\tau ' : G_{\eta}\to G_{\eta}$ is defined by $\tau ' : g\pi\mapsto g\tau\pi$,
and $\tau '$ is the unique $\tau^* : G_{\eta}\to G_{\eta}$ such that $\tau\pi = \pi\tau^*$.

\endproclaim
\demo{Proof} Part (1) follows from the definition of Bruck folders.
 Then (1) and (\cite{A2}, 4.3) imply (2) and (3).

Assume the hypotheses of (4).
 Then $K_{\mu}\subseteq K\subseteq K(\tau)$, so $\tau$ acts on $K_{\mu}$ and
 $H_{\tau}\leq H\leq C_G(\tau)$, so $\tau$ acts on $H_{\mu}$.
 Therefore $\tau$ acts on $G_{\mu}=H_{\mu}K_{\mu}$.
 By (\cite{A2}, 6.2), $\mu$ is a Bol loop folder, and by construction $H_{\mu}\subseteq C_{G_{\mu}}(\tau)$
and $K_{\mu}\subseteq K(\tau)$, so by (\cite{A2}, 5.2),
$\tau_{|\langle K_{\mu}\rangle }=\tau_{\langle K_{\mu}\rangle }$, completing the proof of (4).

Finally assume the hypotheses of (5).
 Then $\xi$ is a normal subfolder of $\xi$, so $\tau$ acts on $G_0$ by (4) and hence induces
$\tau ' : G_{\eta}\to G_{\eta}$ defined by $\tau ' : g\pi\mapsto g\tau\pi$; further $\tau '$ is the
unique map $\tau^* : G_{\eta}\to G_{\eta}$ such that $\tau\pi = \pi\tau^*$.
 As $\pi : K\to K_{\eta}$ is surjective and $K\subseteq K(\tau)$, this implies
$K_{\eta}\subseteq K(\tau ')$, and similarly $\tau '$ centralizes $H_{\eta}$, so by (\cite{A2}, 5.2),
$\Xi_{K_{\eta}}(\langle K_{\eta}\rangle )= 1$ and $\tau ' = \tau_{\eta}$. Thus (5) holds.
\enddemo
\smallpagebreak

\proclaim\nofrills{(5.4)\usualspace } (1) $\xi_{\tau}$ is a Bruck loop folder.

(2) $X_{\tau}=l(\xi_{\tau})$ is of exponent $2$.

(3) $K_{\tau}$ is $G_{\tau}$-invariant and $\Lambda_{\tau}=C_{\Lambda}(\tau)=\tau K_{\tau}$.
\endproclaim
\demo{Proof} By (\cite{A2}, 6.6.5), $\xi_{\tau}$ is a subfolder of $\xi$ and (2) holds.
 Then (1) follows from 5.3.4.
 For $g\in G_{\tau}$ and $k\in K_{\tau}$, $k=\tau\lambda$ for some $\lambda\in\Lambda_{\tau}$,
so $k^g=\tau\lambda^g\in \tau\Lambda_{\tau}=K_{\tau}$ as $\Lambda$ is $G$-invariant, so
$K_{\tau}$ is $G_{\tau}$-invariant. Thus (3) holds.
\enddemo
\smallpagebreak

\proclaim\nofrills{(5.5)\usualspace } Let $U\subseteq H$ and $X=l(\xi)$. Then

(1) $\xi_U$ is a Bruck loop folder.

(2) $X_U=Fix_X(U)$ is a Bruck subloop of $X$ with $l(\xi_U)=X_U$.

(3) $K_U=N_K(U)$.

(4) $\langle K_U\rangle$ is transitive on $X_U$.

(5) $\Lambda_U=N_{\Lambda}(U)=C_{\Lambda}(U)=\tau K_U$.

(6) If $h\in H$ and $h^2\not= 1$ then $h$ is not inverted by any member of $\Lambda$.
 In particular $\tau$ inverts no conjugate of $h$.
\endproclaim
\demo{Proof} By 5.3.1, $\xi$ is an $A_r$-loop folder, so by (\cite{A2}, 4.3.3),
$\xi_U$ is a subfolder of $\xi$. Thus (1) follows from 5.3.4.
 By parts (1) and (2) of (\cite{A2}, 4.3), $X_U$ is a subloop of $X$ with $l(\xi_U)=X_U$.
 Then as $\xi_U$ is a Bruck folder, $X_U$ is a Bruck loop, so (2) holds.
 Parts (3) and (4) follow from parts (4) and (6) of (\cite{A2}, 4.3).

Next $N_{\Lambda}(U)=\tau N_K(U)=\tau K_U$, so (5) follows from (3). Then (6) follows from (5).
\enddemo
\smallpagebreak

\proclaim\nofrills{(5.6)\usualspace } Assume $\xi$ is an envelope, $Q\leq G$ with $Q\ \trianglelefteq\ G^+$,
and set $G^{+*}=G^+/Q$. Then the following are equivalent:

(1) $|k^*|$ is odd for each $k\in K$.

(2) $|\tau^*\lambda^*|$ is odd for each $\lambda\in\Lambda$.

(3) $|G^*|$ is odd.
\endproclaim
\demo{Proof} As $K=\tau\Lambda$, (1) and (2) are equivalent. Trivially (3) implies (2).
 Finally if (2) holds then $\Lambda^* = \tau^{*G}$ and by Glauberman's $Z^*$-Theorem \cite{G1},
$G^*=\langle K^*\rangle = \langle\tau^*\Lambda^*\rangle$ is of odd order.
\enddemo
\smallpagebreak

Recall that a loop $X$ is a $2${\it -loop} if $|X|$ is a power of $2$,
and $X$ is a $2'${\it -loop} if $|X|$ is odd.
\smallpagebreak

\proclaim\nofrills{(5.7)\usualspace } Assume $\xi$ is an envelope and set $X=l(\xi)$. Then

(1) The following are equivalent:

\ \ (a) $X$ is a $2$-loop.

\ \ (b) $G$ is a $2$-group.

\ \ (c) $\alpha\beta$ is a $2$-element for all $\alpha ,\beta\in\Lambda$.

(2) The following are equivalent:

\ \ (a) $X$ is a $2'$-loop.

\ \ (b) $|G|$ is odd.

\ \ (c) $|k|$ is odd for all $k\in K$.
\endproclaim
\demo{Proof} As $|X|=|K|=|G:H|$, (b) implies (a) and (c) in (1) and (2).

Assume $X$ is a $2$-loop. Then $|G:H|$ is a power of $2$, so for each odd prime
$p$ and each element $g$ of order $p$ in $G$, $g$ is conjugate to an element of $H$.
 Thus no member of $\Lambda$ inverts $g$ by 5.5.6.  Hence by the Baer-Suzuki Theorem
(cf. \cite{FGT}, 39.6), $\Lambda\subseteq O_2(G^+)$, so $K=\tau\Lambda\subseteq O_2(G)$.
 Therefore $G=\langle K\rangle$ is a $2$-group.
 Similarly if (1c) holds then no member of $\Lambda$ inverts a nontrivial element of odd order,
so the same argument shows $G$ is a $2$-group, completing the proof of (1).

Assume $X$ is a $2'$-loop. As each $k\in K^{\#}$ is fixed point free on $X$, while
$\langle k\rangle\subseteq K$ by (\cite{A2}, 5.1), $k$ is semiregular on $X$, so $|k|$ divides $|X|$ and hence $|k|$ is odd.
 Therefore (2a) implies (2c), while (2c) implies (2b) by 5.6, completing the proof of (2).
\enddemo
\smallpagebreak

\proclaim\nofrills{(5.8)\usualspace } Let $L$ be a group of odd order and $t$ an involutory automorphism of $L$.
 Then $\mu = (L,C_L(t),K_L(t))$ is a Bruck loop folder, where $K_L(t)=\{l\in L : l^t=l^{-1}\}$.
\endproclaim
\demo{Proof} Let $C=C_L(t)$ and $K=K_L(t)$.
 The map $\sigma : Cg\mapsto [t,g]$ is a well defined injection of $L/C$ into $K$.
 Further as $|L|$ is odd, for $k\in K$, $tk\in t^L$, so $k\in tt^L\subseteq \sigma(L/C)$.
 Thus $\sigma$ is a bijection so $|L:C|=|K|$.
 Finally if $a,b,c\in t^G$ with $1\not= x=ab\in C_G(c)$ then $\langle a\rangle$ and $\langle c\rangle$
are Sylow in the normalizer of $X=\langle x\rangle$, so there is $g\in N_G(X)$ with $a^g=c$ by Sylow's Theorem.
 This is impossible as $a$ inverts $X$, while $c$ centralizes $X$.
 Thus $\mu$ is a Bol loop folder and $\Xi_K(\langle K\rangle )=1$ by the equivalence of parts (1) and (6)
of (\cite{A2}, 6.4). Then by construction, $\mu$ is a Bruck folder.
\enddemo
\smallpagebreak

\proclaim\nofrills{(5.9)\usualspace } Assume $G_0\ \trianglelefteq\ G$, set $G^*=G/G_0$, and
assume $|k^*|$ is odd for each $k\in K$.
 Let $H_0=H\cap G_0$, $K_0=K\cap G_0$, $\xi_0=(G_0,H_0,K_0)$, and $\pi : G\to G^*$ the natural map.
 Then

(1) $G^*$ is of odd order.

(2) $\xi_0$ is a normal subfolder of $\xi$ and $\pi : \xi\to \xi/\xi_0 = (G^*,H^*,K^*)$ is a
surjective morphism of loop folders with $\xi_0=\ker(\pi)$.

(3) $l(\xi)/l(\xi_0)\cong l(\xi/\xi_0)$ is a $2'$-loop.
\endproclaim
\demo{Proof} Part (3) follows from (1), (2), (\cite{A2}, 2.7), and 5.7.2; part (1) follows from 5.6.

Let $t$ be the involutory automorphism of $G^*$ induced by $\tau$ as in 5.3.5, and $K_t=K_{G^*}(t)$.
 By 5.8, $\xi^* = (G^*,C_{G^*}(t),K_t)$ is a Bruck loop folder.
 Further $H^*\leq C_{G^*}(t)$, $K^*\subseteq K_t$, and $G^*=H^*K^*$, so it follows
that $H^*=C_{G^*}(t)$ and $K^*=K_t$. Thus $\xi^*=(G^*,H^*,K^*)$ is a loop folder and
$\pi : \xi\to\xi^*$ is a surjective morphism of folders, so $\xi_0=\ker(\pi)$ is a normal subloop of $\xi$
and $\xi^*=l(\xi)/l(\xi_0)$ by definition of the notation in (\cite{A2}, 2.6). That is (2) holds.
\enddemo
\smallpagebreak

\proclaim\nofrills{(5.10)\usualspace } Assume $|G|$ is odd and let $X=l(\xi)$. Then

(1) $\Lambda = \tau^G$, $H=C_G(\tau)$, and $K=K(\tau)$.

(2) The map $\varphi : J\mapsto (J,C_J(\tau),K_J(\tau))$ is a bijection between the set
$\Cal J$ of $\tau$-invariant subgroups $J$ of $G$ and the set $\Cal F$ of subfolders of $\xi$.

(3) Under the bijection $\varphi$, normal subgroups of $G$ correspond to normal subfolders of $\xi$.

(4) The map $Y\mapsto \langle\kappa(Y)\rangle$ is a bijection between the set of subloops of $X$
and the set $\Cal L$ of $L\in\Cal J$ such that $L=[L,\tau ]$.

(5) $G$ and $X$ are solvable.
\endproclaim
\demo{Proof} The proof of 5.9 in the special case where $G_0=1$ shows that (1) holds.

By 5.8, $\varphi$ is a map from $\Cal J$ into $\Cal F$, and by construction, $\varphi$ is injective.
 If $\mu = (J,H_J,K_J)\in\Cal F$ then $J$ is $\tau$-invariant and $\mu$ is a Bruck folder with
$\tau_{\mu}=\tau_{|J}$ by 5.3.4. Thus $H_J=C_J(\tau)$ and $K_J=K_J(\tau)$ by (1).
 Hence $\varphi$ is a surjection, completing the proof of (2).
 If $J\ \trianglelefteq\ G$, then $\varphi(J)\ \trianglelefteq\ \xi$ by 5.9, so (3) also holds.

Let $\Cal Y$ be the set of subloops of $X$ and $\psi$ the map in (4).
 (cf. Convention 1.9 in \cite{A2} for the definition of $\kappa$.)
 Then $\psi$ is an injection from $\Cal Y$ into $\Cal L$. Further if $L\in\Cal L$ then $\varphi(L)\in\Cal F$,
so $\phi(L)=l(\varphi(L))\in\Cal Y$, and hence $\phi $ is a map from $\Cal L$ to $\Cal Y$.
 Next $\kappa(Y)=K\cap\psi(Y)$, so $\kappa(Y)$ is the set of translations of $\varphi(\psi(Y))$, and hence $Y=l(\varphi(\psi(Y)))=\phi(\psi(Y))$. Similarly $K_L(\tau_{|L})=\kappa(\phi(L))$, so $L=\psi(\phi(L))$,
completing the proof of (4).

By the Odd Order Theorem \cite{FT}, $G$ is solvable. Thus a minimal normal subgroup $L$ of $G$
is an elementary abelian $p$-group for some prime $p$. By (3), $\varphi(L)\ \triangleleft\ \xi$, so by
(\cite{A2}, 2.9), $\phi(L)\ \trianglelefteq\ X$.  As $L$ is abelian, $H\cap L\ \trianglelefteq\ L$, so
$\phi(L)\cong L/H\cap L$ by (\cite{A}, 2.10). By induction on the order of $X$, $X/\phi(X)$ is
solvable, so $X$ is solvable.
\enddemo
\smallpagebreak

\proclaim\nofrills{(5.11)\usualspace } $D(G)=\langle K\rangle\ \trianglelefteq\ G$.
\endproclaim
\demo{Proof} This holds as $G=HK$ and $H$ acts on $K$ via conjugation.
\enddemo
\smallpagebreak

\proclaim\nofrills{(5.12)\usualspace } If $X=O_2(X)\times O(X)$ then $X$ is solvable.
\endproclaim
\demo{Proof} By 5.10.5, $O(X)$ is solvable, while by (\cite{A2}, 7.4) and 5.7.1,
$O_2(X)$ is solvable.
\enddemo
\smallpagebreak

\proclaim\nofrills{(5.13)\usualspace } The following are equivalent:

(1) $X$ is a $2$-element loop.

(2) $k$ is a $2$-element for each $k\in K$.

(3) $\tau\in O_2(G^+)$.
\endproclaim
\demo{Proof} Parts (1) and (2) are equivalent by 4.1.1.
 If (2) holds then $\tau\tau^g$ is a 2-element for each $g\in G$, so (3) holds by the Baer-Suzuki Theorem
 (cf. \cite{FGT}, 39.6). Conversely if (3) holds then for each $\lambda\in\Lambda$,
$\tau\lambda\in \langle\lambda\rangle O_2(G^+)$, so $k=\tau\lambda$ is a $2$-element; that is (3) implies (2).
\enddemo
\smallpagebreak

\subheading{Section 6}
 The proof of Theorem 1
\smallpagebreak

In this section we establish Theorem 1.
 Thus we assume $X$ is a finite Bruck loop and we set $\xi=\epsilon(X)=(G,H,K)$.

Let $X_r$ be the set of $r$-elements of $X$ for $r\in\{2,2'\}$, $K_r=R(X_r)$, and $G_r=\langle K_r\rangle$.
\smallpagebreak

\proclaim\nofrills{(6.1)\usualspace } For each $x\in X$, $x=x_2\circ x_{2'}$ with $x_r\in X_r\cap\langle x\rangle$,
and this expression is unique.

(2) $R(x)=R(x_2)R(x_{2'})$ with $R(x_r)\in K_r$.

(3) $G=G_2 * G_{2'}$ is the central product of $G_2$ and $G_{2'}$;
 that is $G=G_2 G_{2'}$ and $[G_2, G_{2'}]=1$.

(4) $O^{2'}(X)\ \trianglelefteq\ X$.

(5) $G_{2'}\leq O(G)$.

(6) $X_{2'} = O^2(X) = O(X)\ \trianglelefteq\ X$ and $X/O(X)$ is a $2$-element loop.
\endproclaim
\demo{Proof} Parts (1) and (2) follows from 4.1.1 and the corresponding statement for groups.
 Then (3) follows from (1), (2), and Theorem 2.

By (3), $G_2\ \trianglelefteq\ G$. Let $G^*=G/G_2$, $J_r = K\cap G_r$, and
$Y=R^{-1}(J_2)$. For $R(x)\in J_2$, $R(x) = R(y_1)\cdots R(y_n)$ with
each $y_i\in X_2$. Applying both sides to $1\in X$ shows that
$x\in \langle X_2\rangle = O^{2'}(X)$, and so $Y=O^{2'}(X)$.

Each element of $K^*$ is of odd order, so 5.9 tells us that $G^*$ is of odd order and
$Y=O^{2'}(X)$ is a normal subloop of $X$. Let $L=G_{2'}$ and $U=L\cap G_2$.
 Then $L^*\cong L/U$ is of odd order, and hence solvable by the Odd Order Theorem \cite{FT}.
 Also $U\leq Z(L)$ by (3), so $L$ is solvable and hence $L=L_0U$, where $L_0$ is a Hall $2'$-subgroup of $L$
 by Phillip Hall's Theorem 18.5 in \cite{FGT}.
 Then as $U\leq Z(L)$, $L=O^2(L)=L_0$ is of odd order, establishing (5).

By 5.8, $\xi_{2'}=(L,C_L(\tau),K_{2'})$ is a subfolder of $\xi$, so by (\cite{A2}, 1.9),
$X_{2'}=R^{-1}(K_{2'})$ is a subloop of $X$. Then by definition, $X_{2'}=O^2(X)$
 Let $u\in K_{2'}$ and $v\in K$. By (2), $v=ba$ with $b\in K_{2'}$ and $a\in K_2$,
 so $uv=uba$. As $\xi_{2'}$ is a subfolder, $ub=hb'$ with $h\in C_L(\tau)$ and $b'\in K_{2'}$.
 By Theorem 2, $b'a=k'\in K$, so $uv=hv'$; that is the normality condition (NC) of (\cite{A2}, \S2) is satisfied.
 Hence (6) follows from (\cite{A2}, 2.9).
\enddemo
\smallpagebreak

We are now in a position to prove Theorem 1.
 We apply 2.1 to $G_{2'}$ and $G_2$ in the roles of the groups ``$G_1$" and ``$G_2$" in that lemma.
 From the proof of 6.1, the subfolders $\xi_r=(G_r,H\cap G_r,J_r)$ are normal with $J_{2'}=O(X)$
and $J_2=O^{2'}(X)$. By 6.1.3, $G=G_1 G_2=G_3$. Thus by 2.1, $\bar G = \bar G_1\times\bar G_2$
and $X/X_0\cong X_1/X_0\times X_2/X_0=O^{2'}(X)/X_0\times O(X)/X_0$.
 Also $X_0=R^{-1}(J)$, were $J=K_{2'}\cap G_2$.
 By 6.1.3, $J\leq Z(G)$. For $k\in K$, $k=ba$ with $a\in K_2$, $b\in K_{2'}$.
 As $\tau$ inverts $j$ and $b$ and $j\in Z(G)$, $\tau$ inverts $jb$, so $jb\in K_{2'}$ by 5.8.
 Thus $jk=jba\in K$ by Theorem 2.
 Hence $jK\subseteq K$, so $R^{-1}(j)\in Z(X)$; that is $X_0\leq Z(X)$.

We have established the first four statements in Theorem 1, so it remains to establish the fifth.
 Thus we may assume $X$ is solvable. Moreover we assume $X$ is a counterexample of minimal
order to part (5) of Theorem 1. Therefore $O^{2'}(X)\not= O_2(X)$, so $O^{2'}(X)$ is not a $2$-loop.
 Also each proper section $Y$ of $X$ is solvable, so by minimality of $X$, $Y=O_2(Y)\times O(Y)$
and $G_Y=O_2(G_Y)\times O(G_Y)$, where $G_Y$ is the enveloping group of $Y$.
 In particular $X=O^{2'}(X)$ is not a $2$-loop and $G=O^{2'}(G)$ is not a $2$-group.

Suppose $O(X)\not= 1$. By minimality of $X$, $Y=X/O(X)$ is a $2$-loop and $G_Y$ is a $2$-group.
 Let $U=R(O(X))$ and $G^*=G/U$. As $O(X)\leq Z(X)$, $\xi_0=(U,1,U)$ is a normal subfolder of $\xi$, $\xi^*=(G^*,H^*,K^*)=\xi/\xi_0$, and $Y\cong l(\xi/\xi_0)$, so $G_Y=G^*/\ker_{H^*}(G^*)$.
 By 3.3.2, $\ker_{H^*}(G^*)\leq Z(G^*)$, so as $G_Y$ is a $2$-group and $U\leq Z(G)$, $G$ is solvable
and $G=T Z(G)$ for $T\in Syl_2(G)$ by coprime action (cf. \cite{FGT}, 18.7.4).
 Thus $G=O^{2'}(G)=T$, contradicting $G$ not a $2$-group.

Let $Y$ be a maximal normal subloop of $X$. As $O(X)=1$ and $Y=O_2(Y)\times O(Y)$
with $O(Y)=O^2(Y)\leq O^2(X)= O(X)$, $Y$ is a $2$-loop. As $X$ is a solvable $2$-element loop
and $Y$ a maximal normal subloop of $X$, $X/Y\cong\bold Z_2$. Thus $|X|=2|Y|$ is a power of $2$,
so $X$ is a $2$-loop, for our final contradiction.
\smallpagebreak

\subheading{Section 7}
 The proof of Theorem 3.
\smallpagebreak

In this section we establish Theorem 3.

Assume $X$ is an M-loop and let $\xi=\epsilon(X)=(G,H,K)$.
\smallpagebreak

\proclaim\nofrills{(7.1)\usualspace } (1) $X$ is simple.

(2) $O(X)=O_2(X)=1$.

(3) $X$ is a $2$-element loop.

(4) Theorem 3 holds if $X$ is of exponent $2$.
\endproclaim
\demo{Proof} If $Y$ is a proper nontrivial normal subloop of $X$ then $Y$ and $X/Y$
are proper sections of $X$, and hence are solvable; but then $X$ is also solvable, contradicting
the hypothesis that $X$ is an $M$-loop. Therefore (1) holds.
 Then (2) follows from (1) and 5.12, and (3) follows from (2) and Theorem 1.
 Finally if $X$ is of exponent 2 then $X$ is an N-loop, as defined in \cite{A2}, so
Theorem 3 holds in this case by the Main Theorem of \cite{A2}.
\enddemo
\smallpagebreak

By 7.1.3, $X$ is a $2$-element loop. We can repeat many of the lemmas from
(\cite{A2}, \S12), proved there under the stronger hypothesis that $X$ is of exponent $2$.
By 7.1.4, we may assume $X$ is not of exponent $2$.
 Adopt notation 5.2, and for $U\leq G$ set $D(U)=\langle K\cap U\rangle$.
 For $U\leq H$, let $D_U=D(G_U)$.
\smallpagebreak

\proclaim\nofrills{(7.2)\usualspace } Assume $p$ is an odd prime divisor of $|H|$ and let
$1\not= P$ be a $p$-subgroup of $H$. Then

(1) $H$ contains a Sylow $p$-subgroup of $H$.

(2) No member of $\Lambda$ inverts an element of order $p$.

(3) $|N_G(P):N_H(P)|$ is a power of $2$.
\endproclaim
\demo{Proof} As $\xi_P$ is a proper subfolder, $G_P=H_PD_P$,
 and $D_P$ is a $2$-group as $X$ is a $2$-element M-loop.
 Thus (3) holds. Then (3) implies (1), while (1) and 5.5.6 imply (2).
\enddemo
\smallpagebreak

\proclaim\nofrills{(7.3)\usualspace } $|G:H|$ is not a power of $2$.
\endproclaim
\demo{Proof} By 7.1.2, $G$ is not a $2$-group, so as $|X|=|G:H|$, the lemma follows from 5.7.1.
\enddemo
\smallpagebreak

During the remainder of this section we work in the following setup:
\smallpagebreak

 \flushpar {\bf Hypothesis 7.4.} $M^+$ is a maximal overgroup of $\langle\tau\rangle H$ in $G^+$.
 Set $M=M^+\cap G$, $J^+=\ker_{M^+}(G^+)$, $K_M=K\cap M^+$, $\Lambda_M=\tau K_M$,
$D=D(M)$, and $G^{+*}=G^+/J^+$.
\smallpagebreak

\proclaim\nofrills{(7.5)\usualspace } (1) Hypothesis 3.1 is satisfied.

(2) $O_2(G)\langle\tau\rangle = O_2(G^+)\leq J$.

(3) $|G^+:M|$ is even.
\endproclaim
\demo{Proof} Visibly Hypothesis 7.4 implies Hypothesis 3.1, so (1) holds.
 As $X$ is a $2$-element loop, $\tau\in O_2(G^+)$ by 5.13, so (2) holds.
 Finally $\tau\in\Lambda\cap J$ by (2), so (1) and 3.2.5 imply (3).
\enddemo
\smallpagebreak

\proclaim\nofrills{(7.6)\usualspace } (1) $M = HD$ and $D$ is a 2-group.

(2) $|\Lambda_M|=|M:H|$ is a power of $2$.

(3) $|G:M|$ is even but not a power of $2$.

(4) $M$ and $H$ are not $2$-groups.

(5) $D\ \trianglelefteq\ G$, $H^*=M^*$, and $K^*\cap M^*=1$.

(6) Let $N$ be the preimage in $G$ of $F^*(G^*)$. Then $G=HN$.

(7) $J$ is a $2$-group.
\endproclaim
\demo{Proof} The proofs of (1)-(4) are the same as that of the corresponding parts of (\cite{A2}, 12.6).
 Similarly if $D\ \trianglelefteq\ G$ then the proof of (\cite{A2}, 12.6.5) shows that (5) holds, so suppose
$D$ is not normal in $G$. Recall $D\ \trianglelefteq\ M^+$ by 5.11, so $M^+=N_{G^+}(D)$ by maximality
of $M^+$. Next let $D^+=D\langle\tau\rangle$; thus $D=D^+\cap G$ and $D^+=\langle\Lambda_M\rangle$,
with $\Lambda_M=\Lambda\cap M$. Thus $D^+=\langle N_{\Lambda }(D^+)\rangle$,
so taking $D^+\leq T^+\in Syl_2(G^+)$, $D^+=\langle\Lambda\cap N_{T^+}(D^+)\rangle$
and hence $D^+=\langle T^+\cap\Lambda\rangle$, so $T^+\leq N_G(D^+)\leq M^+$, contrary to (3).

Now parts (6) and (7) follows as in the proof of the corresponding parts of (\cite{A2}, 12.6).
\enddemo
\smallpagebreak

\proclaim\nofrills{(7.7)\usualspace } Suppose $1\not= U^*\leq H^*$ is a $p$-group for some odd prime $p$.

(1) $H$ contains a Sylow $p$-group of the preimage of $U^*$ in $G^*$.

(2) $N_{G^*}(U^*)=N_G(P)^*-N_H(P)^*D(C_G(P^*))$.

(3) The triple $G^*$, $M^*$, $K^*$ satisfies Hypothesis N of (\cite{A2}, \S10).

(4) $K^*=\Lambda^*$.
\endproclaim
\demo{Proof} The proofs of (1)-(3) are the same as those of the corresponding parts of (\cite{A2}, 12.7).
 Use 3.3 in proving (3). Note as $\tau\in J$, $K^*=\Lambda^*$ is a union of conjugacy classes of involutions
of $G^*$.
\enddemo
\smallpagebreak

\proclaim\nofrills{(7.8)\usualspace } $F^*(G^*)$ is a nonabelian simple group and $G^*=F^*(G^*)H^*$.
\endproclaim
\demo{Proof} The proof is the same as that of (\cite{A2}, 12.8).
\enddemo
\smallpagebreak

\proclaim{Theorem 7.9} (1) $G^*\cong PGL_2(q)$ with $q=2^n+1$,
$H^*$ is a Borel subgroup of $G^*$, and $K^{*\#}$ consists of the involutions in $G^*-F^*(G^*)$.

(2) $|K^{*\#}|=m=q(q-1)/2$ and $|G:M|=q+1$.

(3) Let $n_0=|K\cap J|$ and $n_1=|K\cap aJ|$ for $a\in K-J$.
 Then $n_0=|M:H|=|D:D\cap H|$ is a power of 2, $n_0=n_12^{n-1}$, and
$|X|=|K|=(q+1)n_0=n_12^n(2^{n-1}+1)$.

(4) $F^*(G)=J=O_2(G)$.
\endproclaim
\demo{Proof} The proofs of (1)-(3) are the same as the corresponding parts of (\cite{A2}, 12.9).
 Note that since $\tau\in J^+$, for $\lambda\in\Lambda - J^+$, $\lambda = \tau a$ for some
$a\in K-J^+$ and the map $k\mapsto \tau k$ is a bijection of $K\cap aJ^+$ with $\lambda J^+\cap\Lambda$.

By 5.4, $G_{\tau}=HD_{\tau}$, where $D_{\tau}=D(G_{\tau})$ and $X_{\tau}$ is of exponent $2$.
 As we are assuming that $X$ is not of exponent $2$, $X_{\tau}\not= X$, so $X_{\tau}$ is a $2$-loop,
and hence $D_{\tau}$ is a $2$-group. Thus each subgroup of $G_{\tau}$ of odd order is fused into $H$
under $D_{\tau}$. However if (4) fails then as $\tau\in J^+=O_2(G^+)$, $G^+=LJ^+$, where $L=E(G)$,
and $\tau$ centralizes $L$. But then $L$ centralizes $D_{\tau}$, so each subgroup of $L$ of odd order
is contained in $H$. Therefore $E=O^2(E)\leq H$, so $E\leq\ker_H(G)=1$, a contradiction.
\enddemo
\smallpagebreak

Observe that 7.1 and Theorem 7.9 establish Theorem 3.
\smallpagebreak

\subheading{Section 8}
 The proof of Corollary 4
\smallpagebreak

In this section we prove Corollary 4.
 Thus we assume $X$ is a Bol loop which is also an $A_r$-loop.
 Let $\xi =\epsilon(X)=(G,H,K)$, $G^*=G/\Xi_K(G)$, and $\hat G= G^*/\ker_{H^*}(G^*)$.
\smallpagebreak

\proclaim\nofrills{(8.1)\usualspace } (1) $\xi_{\Xi}=(\Xi_K(G),1,\Xi_K(G))$ is a normal subfolder of $G$.

(2) $\Xi(X)$ is a normal subloop of $X$.

(3) $\Xi(X)$ is isomorphic to the group $\Xi_K(G)$.

(4) $X/\Xi(X)\cong l(G^*,H^*,K^*)$.

(5) There is a unique automorphism $\tau$ of $G^*$ with $\tau^2=1$ and $K^*\subseteq K(\tau)$.

(6) $H^*\leq C_{G^*}(\tau )$.

(7) $\epsilon(X/\Xi(X))\cong (\hat G,\hat H,\hat K)$.

(8) $X/\Xi(X)$ is a Bruck loop.

(9) $\ker_{H^*}(G^*)\leq Z(G^*)$.
\endproclaim
\demo{Proof} Parts (1)-(4) are the corresponding parts of (\cite{A2}, 6.5).
 Then (7) is a consequence of (4) and (\cite{A2}, 2.9.2).
 Part (5) follows from (\cite{A2}, 5.1.3).
 As $X$ is an $A_r$-loop, $H$ acts on $K$ via conjugation by
(\cite{A2}, 4.1), so $H^*$ acts on $K^*$.
 Therefore (6) follows from the uniqueness of $\tau$ in (5).
 Thus $\tau$ acts on $\ker_{H^*}(G^*)$ and hence induces an automorphism
$\hat\tau$ of $\hat G$ centralizing $\hat H$ with $\hat K\subseteq K(\hat\tau )$.
 Therefore $\Xi_{\hat K}(\hat G)=1$ by (\cite{A2}, 5.1.3).
 Then (8) follows from this fact and the fact that $\hat\tau$ centralizes $\hat H$,
given the equivalence of parts (1) and (4) of 5.1. Finally (9) follows from (\cite{A2}, 4.3.4).
\enddemo
\smallpagebreak

Assume $X$ is solvable. Then $\Xi(X)$ is solvable, so the subgroup $\Xi_K(G)$
is also solvable by 8.1.3. Next by 8.1.8, $X/\Xi(X)$ is a Bruck loop, and solvable
as $X$ is solvable. Then by Theorem 1 and the Odd Order Theorem \cite{FT},
the enveloping group $\bar G$ of $X/\Xi(X)$ is solvable, and by 8.1.7, $\bar G\cong\hat G$.
 Then appealing to 8.1.9, $G$ is solvable.

Thus to complete the proof of Corollary 3 it remains to show that $X$ is solvable if $G$ is solvable.
 Assume otherwise and choose a counter example $X$ of minimal order.

As $G$ is solvable, so are $\Xi_K(G)$ and $\hat G$. Hence $\Xi(X)$ is solvable by 8.1.3.
 Further by 8.1.8, $X/\Xi(X)$ is a Bruck loop, and by 8.1.7, $\hat G$ is its enveloping group.
 Thus if $\Xi(X)\not= 1$ then $X/\Xi(X)$ is solvable by minimality of $X$, so $X$ is solvable,
contrary to the choice of $X$. Therefore $X$ is radical free and a Bruck loop.

Next if $Y$ is a proper section of $X$, then by (\cite{A2}, 2.9), the enveloping group of $Y$
is a section of $G$, and hence is solvable.  Therefore by minimality of $X$, $Y$ is solvable.
 Thus $X$ is an M-loop.  But now Theorem 3 supplies a contradiction, since $G$ is solvable.
\smallpagebreak

\enddocument

\Refs
\widestnumber\key{FGT}

\ref
 \key{A1}
 \by M. Aschbacher
 \paper Near subgroups of finite groups
 \jour J. Group Theory \vol 1 \yr 1998 \pages 113--129
\endref
\ref
 \key{A2}
 \by M. Aschbacher
 \paper On Bol loops of of exponent 2\toappear
\endref
\ref
 \key {Ba}
 \by R. Baer
 \paper Nets and groups
 \jour Trans. Amer. Math. Soc. \vol 47 \yr 1939 \pages 110--141
\endref
\ref
 \key{Br}
 \by R. Bruck
 \book A Survey of Binary Systems
 \publ Springer-Verlag \publaddr Berlin \yr 1971
\endref
\ref
 \key{FT}
 \by W. Feit and J. Thompson
 \paper Solvability of groups of odd order
 \jour Pacific J. Math. \vol 13 \yr 1963 \pages 755--1029
\endref
\ref
 \key{FGT}
 \by M. Aschbacher
 \book Finite Group Theory
 \publ Cambridge Univ. Press \publaddr Cambridge \yr 1986
\endref
\ref
 \key{FU}
 \by T. Foguel and A. A. Ungar
 \paper Gyrogroups and the decomposition of groups into twisted subgroups and subgroups
 \jour  Pacific J. Math. \vol 197  \yr 2001 \pages 1--11
\endref
\ref
 \key{G1}
 \by G. Glauberman
 \paper Central elements in core free groups
 \jour J. Algebra \vol 4 \yr 1966 \pages 403--420
\endref
\ref
 \key{G2}
 \by G. Glauberman
 \paper On loops of odd order, I
 \jour J. Algebra \vol 1 \yr 1964 \pages 374--396
\endref
\ref
 \key{G3}
 \by G. Glauberman
 \paper On loops of odd order, II
 \jour J. Algebra \vol 8 \yr 1968 \pages 393--414
\endref
\ref
 \key{Ki}
 \by H. Kiechle
 \book Theory of K-loops
 \publ Springer Verlag \publaddr Berlin \yr 2002
\endref
\ref
 \key{Pf}
 \by H. O. Pflugfelder
 \book Quasigroups and Loops: Introduction
 \publ Heldermann Verlag \publaddr Berlin \yr 1990
\endref
\ref
 \key{Ro}
 \by D. Robinson
 \paper Bol loops
 \jour Trans. Amer. Math. Soc. \vol 123 \yr 1966 \pages 341--354
\endref
\endRefs
\bye